\documentclass[11pt]{article}
\usepackage[]{amsmath,amssymb}

\newtheorem{theorem}{Theorem}[section]
\newtheorem{lemma}[theorem]{Lemma}
\newtheorem{proposition}[theorem]{Proposition}
\newtheorem{corollary}[theorem]{Corollary}
\newtheorem{exAux}[theorem]{Example}

\newtheorem{Def}[theorem]{Definition}
\newenvironment{definition}{\begin{Def} \rm}{\end{Def}}
\newtheorem{Note}[theorem]{Note}
\newenvironment{note}{\begin{Note} \rm}{\end{Note}}
\newtheorem{Rem}[theorem]{Remark}
\newenvironment{remark}{\begin{Rem} \rm}{\end{Rem}}
\newtheorem{Ass}[theorem]{Assumption}

\newenvironment{proof}{\medskip\noindent{\bf Proof.\ }}{\qed\medskip}
\newenvironment{proofof}[1]{\medskip\noindent{\bf Proof  of {#1}.\ 
}}{\qed\medskip}
\newcommand{\qed}{\hfill\mbox{$\Box$\qquad\qquad}}

\newcommand{\Mat}[1]{\text{\rm Mat}_{#1}(\mathbb{K})}

\begin{document}
\thispagestyle{empty}

\begin{center}
\LARGE \bf
\noindent
Some trace formulae involving the split\\
sequences of a Leonard pair
\end{center}

\smallskip

\begin{center}
\Large
Kazumasa Nomura and Paul Terwilliger
\end{center}

\smallskip

\begin{quote}
\small 
\begin{center}
\bf Abstract
\end{center}
Let $\mathbb{K}$ denote a field, and let $V$ denote a vector
space over $\mathbb{K}$ with finite positive dimension.
We consider a pair of linear transformations
$A:V \to V$ and $A^*:V \to V$ that satisfy (i), (ii) below:
\begin{itemize}
\item[(i)] There exists a basis for $V$ with respect to which the
matrix representing $A$ is irreducible tridiagonal and the matrix
representing $A^*$ is diagonal.
\item[(ii)] There exists a basis for $V$ with respect to which the
matrix representing $A^*$ is irreducible tridiagonal and the matrix
representing $A$ is diagonal.
\end{itemize}
We call such a pair a {\em Leonard pair} on $V$.
Let $\text{diag} (\theta_0, \theta_1,\ldots, \theta_d)$ denote the
diagonal matrix referred to in (ii) above and 
let $\text{diag} (\theta^*_0, \theta^*_1,\ldots, \theta^*_d)$
denote the diagonal matrix referred to in (i) above.
It is known that there exists a basis $u_0, u_1, \ldots, u_d$
for $V$ and there exist scalars 
$\varphi_1, \varphi_2, \ldots, \varphi_d$ in $\mathbb{K}$ such that
$A u_i= \theta_i u_i + u_{i+1}$ $(0 \leq i \leq d-1)$,
$A u_d=\theta_d u_d$, 
$A^* u_i = \varphi_{i}u_{i-1}+\theta^*_i u_i$ $(1 \leq i \leq d)$,
$A^* u_0 = \theta^*_0 u_0$.
The sequence $\varphi_1, \varphi_2, \ldots, \varphi_d$ is called
the {\em first split sequence} of the Leonard pair.
It is known that there exists a basis $v_0, v_1, \ldots, v_d$
for $V$ and there exist scalars $\phi_1, \phi_2, \ldots, \phi_d$ 
in $\mathbb{K}$ such that
$A v_i= \theta_{d-i} v_i + v_{i+1}$ $(0 \leq i \leq d-1)$,
$A v_d=\theta_0 v_d$, 
$A^* v_i = \phi_{i}v_{i-1}+\theta^*_i v_i$ $(1 \leq i \leq d)$,
$A^* v_0 = \theta^*_0 v_0$.
The sequence $\phi_1, \phi_2, \ldots, \phi_d$ is called the
{\em second split sequence} of the Leonard pair.
We display some attractive formulae for the first and second split
sequence that involve the trace function.
\end{quote}

\section{Leonard pairs and Leonard systems}

Throughout this paper $\mathbb{K}$ will denote an arbitrary field.
We begin by recalling the notion of a Leonard pair.
We will use the following notation.
A square matrix $X$ is called {\em tridiagonal}
whenever each nonzero entry lies on either the diagonal, the subdiagonal,
or the superdiagonal. Assume $X$ is tridiagonal.
Then $X$ is called {\em irreducible}
whenever each entry on the subdiagonal is nonzero and each entry on
the superdiagonal is nonzero.

\medskip

\begin{definition}  \cite{T:Leonard}  \label{def:LP}
Let $V$ denote a vector space over $\mathbb{K}$ with finite positive
dimension.
By a {\em Leonard pair} on $V$ we mean an ordered pair of linear transformations
$A:V \to V$ and $A^*:V \to V$ that satisfy (i), (ii) below:
\begin{itemize}
\item[(i)] There exists a basis for $V$ with respect to which the
matrix representing $A$ is irreducible tridiagonal and the matrix
representing $A^*$ is diagonal.
\item[(ii)] There exists a basis for $V$ with respect to which the
matrix representing $A^*$ is irreducible tridiagonal and the matrix
representing $A$ is diagonal.
\end{itemize}
\end{definition}

\begin{note}
It is a common notational convention to use $A^*$ to represent the
conjugate-transpose of $A$. We are {\em not} using this convention.
In a Leonard pair $A$, $A^*$ the linear transformations $A$ and
$A^*$ are arbitrary subject to (i) and (ii) above.
\end{note}

\medskip
We refer the reader to \cite{AC}, \cite{AC2},
\cite{ITT}, \cite{IT:shape},
\cite{IT:uqsl2hat}, \cite{N:aw}, \cite{N:refine}, \cite{N:height1},
\cite{NT:balanced}, \cite{P}, \cite{T:sub1}, \cite{T:sub3}, \cite{T:Leonard},
\cite{T:qSerre}, \cite{T:24points}, \cite{T:conform}, \cite{T:intro},
\cite{T:intro2}, \cite{T:split}, \cite{T:array}, \cite{T:qRacah},
\cite{T:survey}, \cite{TV}, \cite{V}
for background on Leonard pairs. We especially recommend the
survey \cite{T:survey}.

\medskip
When working with a Leonard pair, it is convenient to consider a closely
related object called a {\em Leonard system}. 
To prepare for our definition
of a Leonard system, we recall a few concepts from linear algebra.
Let $d$ denote a nonnegative integer and let
$\Mat{d+1}$ denote the $\mathbb{K}$-algebra consisting of all $d+1$ by
$d+1$ matrices that have entries in $\mathbb{K}$. 
We index the rows and  columns by $0, 1, \ldots, d$. 
For the rest of this paper we let $\cal A$ denote a $\mathbb{K}$-algebra 
isomorphic to $\Mat{d+1}$. 
Let $V$ denote a simple $\cal A$-module. We remark that $V$ is unique
up to isomorphism of $\cal A$-modules, and that $V$ has dimension $d+1$.
Let $v_0, v_1, \ldots, v_d$ denote a basis for $V$.
For $X \in {\cal A}$ and $Y \in \Mat{d+1}$, we say 
{\em $Y$ represents $X$ with respect to} $v_0, v_1, \ldots, v_d$
whenever $X v_j = \sum_{i=0}^d Y_{ij}v_i$ for $0 \leq j \leq d$.
Let $A$ denote an element of $\cal A$. We say $A$ is {\em multiplicity-free}
whenever it has $d+1$ mutually distinct eigenvalues in $\mathbb{K}$. 
Let $A$ denote a multiplicity-free element in $\cal A$.
Let $\theta_0, \theta_1, \ldots, \theta_d$ denote an ordering 
of the eigenvalues of $A$, and for $0 \leq i \leq d$ put
\[
    E_i = \prod_{\stackrel{0 \leq j \leq d}{j\neq i}}
             \frac{A-\theta_j I}{\theta_i - \theta_j},
\]
where $I$ denotes the identity of $\cal A$. 
 We observe
(i) $AE_i = \theta_i E_i$ $(0 \leq i \leq d)$;
(ii) $E_i E_j = \delta_{i,j} E_i$ $(0 \leq i,j \leq d)$;
(iii) $\sum_{i=0}^{d} E_i = I$;
(iv) $A = \sum_{i=0}^{d} \theta_i E_i$.
We call $E_i$ the {\em primitive idempotent} of $A$ associated with
$\theta_i$. 
We now define a Leonard system.

\medskip

\begin{definition}  \cite{T:Leonard}   \label{def:LS}
By a {\em Leonard system} in $\cal A$ we mean a sequence
\[
\Phi=(A; \{E_i\}_{i=0}^d; A^*; \{E^*_i\}_{i=0}^d)
\]
that satisfies (i)--(v) below.
\begin{itemize}
\item[(i)] Each of $A$, $A^*$ is a multiplicity-free element in $\cal A$.
\item[(ii)] $E_0, E_1, \ldots, E_d$ is an ordering of the
   primitive idempotents of $A$.
\item[(iii)] $E^*_0, E^*_1, \ldots, E^*_d$ is an ordering of the
   primitive idempotents of $A^*$.
\item[(iv)] For $0 \leq i,j \leq d$, 
\[
   E_i A^* E_j =
    \begin{cases}  
        0 & \text{\rm if $|i-j|>1$},  \\
        \neq 0 & \text{\rm if $|i-j|=1$}.
    \end{cases}
\]
\item[(v)] For $0 \leq i,j \leq d$, 
\[
   E^*_i A E^*_j =
    \begin{cases}  
        0 & \text{\rm if $|i-j|>1$},  \\
        \neq 0 & \text{\rm if $|i-j|=1$}.
    \end{cases}
\]
\end{itemize}
We say $\Phi$ is {\em over} $\mathbb{K}$.
We call $\cal A$ the {\em ambient algebra} of $\Phi$.
\end{definition}

\medskip

We now recall the notion of isomorphism for Leonard systems.
Let $\Phi$ denote the Leonard system from Definition \ref{def:LS}
and let $\sigma : {\cal A} \to {\cal A}'$ denote an isomorphism
of  $\mathbb{K}$-algebras. We write
$\Phi^{\sigma} :=(A^{\sigma}; \{E^{\sigma}_i\}_{i=0}^d; A^{*\sigma};
     \{E^{*\sigma}_i\}_{i=0}^d)$
and observe $\Phi^{\sigma}$ is a Leonard system in ${\cal A}'$.

\medskip

\begin{definition}
Let $\Phi$ and $\Phi'$ denote Leonard systems over $\mathbb{K}$.
By an {\em isomorphism of Leonard systems from $\Phi$ to $\Phi'$}
we mean an isomorphism of $\mathbb{K}$-algebras $\sigma$ from the
ambient algebra of $\Phi$ to the ambient algebra of $\Phi'$ such that
$\Phi'=\Phi^{\sigma}$. The Leonard systems $\Phi$, $\Phi'$ are said
to be {\em isomorphic} whenever there exists an isomorphism of
Leonard systems from $\Phi$ to $\Phi'$.
\end{definition}

\section{The $D_4$ action}

Let $\Phi=(A; \{E_i\}_{i=0}^d; A^*; \{E^*_i\}_{i=0}^d)$ denote
a Leonard system in $\cal A$.
Then each of the following is a Leonard system in $\cal A$:
\begin{eqnarray*}
\Phi^{*}  &:=& 
       (A^*; \{E^*_i\}_{i=0}^d; A; \{E_i\}_{i=0}^d),  \\
\Phi^{\downarrow} &:=&
       (A; \{E_i\}_{i=0}^d; A^*; \{E^*_{d-i}\}_{i=0}^d), \\
\Phi^{\Downarrow} &:=&
       (A; \{E_{d-i}\}_{i=0}^d; A^*; \{E^*_{i}\}_{i=0}^d).
\end{eqnarray*}
Viewing $*$, $\downarrow$, $\Downarrow$ as permutations on the set of
all Leonard systems in $\cal A$,
\begin{equation}    \label{eq:relation1}
*^2 = \downarrow^2 = \Downarrow^2 = 1,
\end{equation}
\begin{equation}    \label{eq:relation2}
\Downarrow * = * \downarrow, \quad
\downarrow * = * \Downarrow, \quad
\downarrow \Downarrow = \Downarrow \downarrow.
\end{equation}
The group generated by symbols $*$, $\downarrow$, $\Downarrow$ subject
to the relations (\ref{eq:relation1}) and (\ref{eq:relation2}) is the
dihedral group $D_4$. We recall $D_4$ is the group of symmetries of a
square, and has $8$ elements.
Apparently $*$, $\downarrow$, $\Downarrow$ induce an action of $D_4$
on the set of all Leonard systems in $\cal A$.
Two Leonard systems will be called {\em relatives} whenever they are
in the same orbit of this $D_4$ action. Assuming $d \geq 1$ to avoid
trivialities, the relatives of $\Phi$ are as follows:

\medskip
\noindent
\begin{center}
\begin{tabular}{c|c}
name  &  relative \\
\hline
$\Phi$ & 
       $(A; \{E_i\}_{i=0}^d; A^*;  \{E^*_i\}_{i=0}^d)$ \\ 
$\Phi^{\downarrow}$ &
       $(A; \{E_i\}_{i=0}^d; A^*;  \{E^*_{d-i}\}_{i=0}^d)$ \\ 
$\Phi^{\Downarrow}$ &
       $(A; \{E_{d-i}\}_{i=0}^d; A^*;  \{E^*_i\}_{i=0}^d)$ \\ 
$\Phi^{\downarrow \Downarrow}$ &
       $(A; \{E_{d-i}\}_{i=0}^d; A^*;  \{E^*_{d-i}\}_{i=0}^d)$ \\ 
$\Phi^{*}$  & 
       $(A^*; \{E^*_i\}_{i=0}^d; A;  \{E_i\}_{i=0}^d)$ \\ 
$\Phi^{\downarrow *}$ &
       $(A^*; \{E^*_{d-i}\}_{i=0}^d; A;  \{E_i\}_{i=0}^d)$ \\ 
$\Phi^{\Downarrow *}$ &
       $(A^*; \{E^*_i\}_{i=0}^d; A;  \{E_{d-i}\}_{i=0}^d)$ \\ 
$\Phi^{\downarrow \Downarrow *}$ &
       $(A^*; \{E^*_{d-i}\}_{i=0}^d; A;  \{E_{d-i}\}_{i=0}^d)$
\end{tabular}
\end{center}

\section{The parameter array}

Let $\Phi=(A; \{E_i\}_{i=0}^d; A^*; \{E^*_i\}_{i=0}^d)$ denote
a Leonard system in $\cal A$.
We now associate with $\Phi$ some parameters.
For $0 \leq i \leq d$ we let $\theta_i$ 
denote the eigenvalue of $A$ associated with
$E_i$, and let  $\theta^*_i$ denote the eigenvalue
of $A^*$ associated with $E^*_i$.
We call $\theta_0, \theta_1, \ldots, \theta_d$ 
 the {\em eigenvalue sequence} of $\Phi$, and
 $\theta^*_0, \theta^*_1, \ldots, \theta^*_d$
the {\em dual eigenvalue sequence} of $\Phi$.
Let $V$ denote a simple $\cal A$-module. 
Recall $V$ has dimension $d+1$.
By \cite[Theorem 3.2]{T:Leonard}  there exists a
basis $u_0, u_1, \ldots, u_d$ for $V$ and there exist nonzero scalars
$\varphi_1, \varphi_2, \ldots, \varphi_d$
in $\mathbb{K}$ such that
\[
   A u_i = \theta_i u_i + u_{i+1}  \quad  (0 \leq i \leq d-1),
   \qquad  Au_d = \theta_d u_d,
\]
\[
   A^* u_i = \varphi_{i}u_{i-1}+\theta^*_i u_i \quad (1 \leq i \leq d),
   \qquad  A^* u_0 =  \theta^*_0 u_0.
\]
The basis $u_0, u_1, \ldots, u_d$ is called a {\em $\Phi$-split basis} for $V$.
The sequence $\varphi_1, \varphi_2, \ldots, \varphi_d$ is uniquely
determined by $\Phi$; 
we call this sequence the {\em first split sequence} for $\Phi$.
The first split sequence for $\Phi^{\Downarrow}$ is denoted by
$\phi_1, \phi_2, \ldots, \phi_d$ and called the
 {\em second split sequence} for $\Phi$.

\medskip

\begin{definition}
Let $\Phi$ denote a Leonard system.
By the {\em parameter array} of $\Phi$ we mean the sequence
\[
      p= (\{\theta_i\}_{i=0}^d; \{\theta^*_i\}_{i=0}^d;
        \{\varphi_i\}_{i=1}^d; \{\phi_i\}_{i=1}^d),
\]
where $\theta_0, \theta_1, \ldots, \theta_d$ denotes the eigenvalue
sequence of $\Phi$, 
$\theta^*_0, \theta^*_1, \ldots, \theta^*_d$ denotes
the dual eigenvalue sequence of $\Phi$, 
$\varphi_1, \varphi_2, \ldots, \varphi_d$ denotes 
the first split sequence of $\Phi$, 
and $\phi_1, \phi_2, \ldots, \phi_d$ denotes 
the second split sequence of $\Phi$.
\end{definition}

\medskip

\begin{theorem}   \cite[Theorem 1.9]{T:Leonard}     \label{thm:classify}
Let 
\[
 p = (\{\theta_i\}_{i=0}^d; \{\theta^*_i\}_{i=0}^d;
        \{\varphi_i\}_{i=1}^d; \{\phi_i\}_{i=1}^d)
\]
denote a sequence of scalars taken from $\mathbb{K}$. 
Then there exists a Leonard system $\Phi$ over $\mathbb{K}$ with
parameter array $p$ if and only if (i)--(v) hold below.
\begin{itemize}
\item[(i)]  $\varphi_i \neq 0$, $\phi_i \neq 0$ $(1 \leq i \leq d)$.
\item[(ii)] $\theta_i \neq \theta_j$, $\theta^*_i \neq \theta^*_j$
   if $i \neq j$ $(0 \leq i,j \leq d$).
\item[(iii)] For $1 \leq i \leq d$,
\[
 \varphi_i = \phi_1 
 \sum_{h=0}^{i-1} \frac{\theta_h - \theta_{d-h}}
                       {\theta_0 - \theta_d}
         + (\theta^*_{i}-\theta^*_{0})(\theta_{i-1}-\theta_{d}).
\]
\item[(iv)] For $1 \leq i \leq d$,
\[
\phi_i = \varphi_1
\sum_{h=0}^{i-1} \frac{\theta_h - \theta_{d-h}}
                       {\theta_0 - \theta_d}
         + (\theta^*_{i}-\theta^*_{0})(\theta_{d-i+1}-\theta_{0}).
\]
\item[(v)] The expressions
\begin{equation}        \label{eq:beta}
   \frac{\theta_{i-2}-\theta_{i+1}}{\theta_{i-1}-\theta_{i}},
 \quad
   \frac{\theta^*_{i-2}-\theta^*_{i+1}}{\theta^*_{i-1}-\theta^*_{i}}
\end{equation}
are equal and independent of $i$ for $2 \leq i \leq d-1$.
\end{itemize}
Suppose (i)--(v) hold. Then $\Phi$ is unique up to isomorphism of
Leonard systems.
\end{theorem}

\medskip

\begin{theorem}   \cite[Theorem 1.11]{T:Leonard}        \label{thm:rel}
Let $\Phi$ denote a Leonard system and let
\[
   p = (\{\theta_i\}_{i=0}^d; \{\theta^*_i\}_{i=0}^d;
        \{\varphi_i\}_{i=1}^d; \{\phi_i\}_{i=1}^d)
\]
denote the parameter array of $\Phi$. Then the following (i)--(iii) hold.
\begin{itemize}
\item[(i)] The paramater array of $\Phi^*$ is $p^*$ where
\[
   p^* = (\{\theta^*_i\}_{i=0}^d; \{\theta_i\}_{i=0}^d;
        \{\varphi_i\}_{i=1}^d; \{\phi_{d-i+1}\}_{i=1}^d).
\]
\item[(ii)]
The parameter array of $\Phi^{\downarrow}$ is $p^{\downarrow}$ where
\[
   p^{\downarrow}
       = (\{\theta_i\}_{i=0}^d; \{\theta^*_{d-i}\}_{i=0}^d;
        \{\phi_{d-i+1}\}_{i=1}^d; \{\varphi_{d-i+1}\}_{i=1}^d).
\]
\item[(iii)]
The parameter array of $\Phi^{\Downarrow}$ is $p^{\Downarrow}$ where
\[
   p^{\Downarrow}
      = (\{\theta_{d-i}\}_{i=0}^d; \{\theta^*_{i}\}_{i=0}^d;
        \{\phi_{i}\}_{i=1}^d; \{\varphi_{i}\}_{i=1}^d).
\]
\end{itemize}
\end{theorem}

\section{Trace formulae involving the split sequences}

In this section we obtain some trace formulae that involve
the split sequences. We use the following notation.
Throughout this section $\Phi=(A; \{E_i\}_{i=0}^d; A^*; \{E^*_i\}_{i=0}^d)$ 
denotes a Leonard system in $\cal A$ and 
$p=(\{\theta_i\}_{i=0}^d; \{\theta^*_i\}_{i=0}^d;
        \{\varphi_i\}_{i=1}^d; \{\phi_i\}_{i=1}^d)$
denotes the corresponding parameter array.
Let $\lambda$ denote an indeterminate and let 
$\mathbb{K}[\lambda]$ denote the $\mathbb{K}$-algebra consisting of all polynomials 
in $\lambda$ that have coefficients in $\mathbb{K}$.

\medskip

\begin{definition}
For $0 \leq i \leq d$ we let $\tau_i$, $\eta_i$, $\tau^*_i$,
$\eta^*_i$ denote the following polynomials in $\mathbb{K}[\lambda]$:
\begin{eqnarray}
 \tau_i &=& (\lambda-\theta_0)(\lambda-\theta_1) \cdots(\lambda-\theta_{i-1}),
      \label{eq:deftau} \\
 \eta_i &=& 
   (\lambda-\theta_d)(\lambda-\theta_{d-1})\cdots (\lambda-\theta_{d-i+1}), 
       \label{eq:defeta} \\
 \tau^*_i &=& (\lambda-\theta^*_0)(\lambda-\theta^*_1) \cdots(\lambda-\theta^*_{i-1}),
    \label{eq:deftaus} \\
 \eta^*_i &=& 
   (\lambda-\theta^*_d)(\lambda-\theta^*_{d-1})\cdots (\lambda-\theta^*_{d-i+1}).
      \label{eq:defetas}
\end{eqnarray}
\end{definition}

\medskip

\begin{lemma}         \label{lem:tauAu0}
Let $V$ denote a simple $\cal A$-module and let $u_0, u_1,\ldots, u_d$
denote a $\Phi$-split basis for $V$. Then
\[
     \tau_i(A) u_0 = u_i  \qquad (0 \leq i \leq d).
\]
\end{lemma}

\begin{proof}
Immediate from (\ref{eq:deftau}) and since
$u_i=(A-\theta_{i-1}I)u_{i-1}$ for $1 \leq i \leq d$.
\end{proof}

\medskip

\begin{lemma}   \cite[Lemma 4.6]{T:Leonard}     \label{lem:entryBs}
For $0 \leq r \leq d$ let $B^*_r$ denote the matrix in $\Mat{d+1}$
that represents $E^*_r$ with respect to a $\Phi$-split basis.
Then $B^*_r$ has $(j,i)$-entry
\[
  \frac{\varphi_1 \varphi_2\cdots \varphi_i
            \tau^*_j(\theta^*_r) \eta^*_{d-i}(\theta^*_r)}
       {\varphi_1 \varphi_2 \cdots \varphi_j
            \tau^*_r(\theta^*_r) \eta^*_{d-r}(\theta^*_r)}
\]
for $0 \leq i,j \leq d$.
\end{lemma}

\medskip

\begin{proposition}        \label{prop:trace}
For $0 \leq i \leq d$ we have
\begin{eqnarray}     
\text{\rm tr}(\tau_{i}(A)E^*_0)
 &=& \frac{\varphi_{1} \varphi_{2} \cdots \varphi_{i}}
     {(\theta^*_{0}-\theta^*_{1})(\theta^*_{0}-\theta^*_{2}) 
                      \cdots (\theta^*_{0}-\theta^*_{i})},
     \label{eq:trace1}   \\
\text{\rm tr}(\tau_i(A)E^*_d)
 &=& \frac{\phi_{d} \phi_{d-1} \cdots \phi_{d-i+1}}
     {(\theta^*_{d}-\theta^*_{d-1})(\theta^*_{d}-\theta^*_{d-2}) 
             \cdots (\theta^*_{d}-\theta^*_{d-i})},
     \label{eq:trace2}   \\
\text{\rm tr}(\eta_{i}(A)E^*_0)
 &=& \frac{\phi_{1} \phi_{2} \cdots \phi_{i}}
     {(\theta^*_{0}-\theta^*_{1})(\theta^*_{0}-\theta^*_{2}) 
                      \cdots (\theta^*_{0}-\theta^*_{i})},
     \label{eq:trace3}   \\
\text{\rm tr}(\eta_i(A)E^*_d)
 &=& \frac{\varphi_{d} \varphi_{d-1} \cdots \varphi_{d-i+1}}
     {(\theta^*_{d}-\theta^*_{d-1})(\theta^*_{d}-\theta^*_{d-2}) 
             \cdots (\theta^*_{d}-\theta^*_{d-i})},
     \label{eq:trace4}   \\
\text{\rm tr}(\tau^*_{i}(A^*)E_0)
 &=& \frac{\varphi_{1} \varphi_{2} \cdots \varphi_{i}}
     {(\theta_{0}-\theta_{1})(\theta_{0}-\theta_{2}) 
                      \cdots (\theta_{0}-\theta_{i})},
     \label{eq:trace5}   \\
\text{\rm tr}(\tau^*_{i}(A^*)E_d)
 &=& \frac{\phi_{1} \phi_{2} \cdots \phi_{i}}
     {(\theta_{d}-\theta_{d-1})(\theta_{d}-\theta_{d-2}) 
                      \cdots (\theta_{d}-\theta_{d-i})},
     \label{eq:trace6}   \\
\text{\rm tr}(\eta^*_{i}(A^*)E_0)
 &=& \frac{\phi_{d} \phi_{d-1} \cdots \phi_{d-i+1}}
     {(\theta_{0}-\theta_{1})(\theta_{0}-\theta_{2}) 
                      \cdots (\theta_{0}-\theta_{i})},
     \label{eq:trace7}   \\
\text{\rm tr}(\eta^*_{i}(A^*)E_d)
 &=& \frac{\varphi_{d} \varphi_{d-1} \cdots \varphi_{d-i+1}}
     {(\theta_{d}-\theta_{d-1})(\theta_{d}-\theta_{d-2}) 
                      \cdots (\theta_{d}-\theta_{d-i})},
     \label{eq:trace8} 
\end{eqnarray}
where tr means trace.
\end{proposition}

\begin{proof}
We first show (\ref{eq:trace1}).
Let $V$ denote a simple $\cal A$-module and let $u_0, u_1, \ldots, u_d$
denote a $\Phi$-split basis for $V$. Let the matrix $B^*_0$ be as in
Lemma \ref{lem:entryBs}.
For $0 \leq i \leq d$ let $F_i$ denote the matrix in $\Mat{d+1}$
that represents $\tau_i(A)$ with respect to the basis $u_0, u_1, \ldots, u_d$.
Recall $B^*_0$ is the matrix in $\Mat{d+1}$ that represents $E^*_0$
with respect to $u_0, u_1, \ldots, u_d$. 
Then $F_i B^*_0$ represents $\tau_i(A)E^*_0$ with respect to
$u_0, u_1, \ldots, u_d$. Therefore $F_iB^*_0$ and $\tau_i(A)E^*_0$
have the same trace. We now compute the trace of $F_iB^*_0$.
First observe by Lemma \ref{lem:tauAu0} that in column $0$ of $F_i$ the
$i$th entry is $1$ and all other entries are $0$.
Next observe by Lemma \ref{lem:entryBs} that in $B^*_0$ the rows
$1, \ldots, d$ are all zero. By the above two observations the trace
of $F_iB^*_0$ is equal to the $(0,i)$-entry of $B^*_0$.
By Lemma \ref{lem:entryBs} the $(0,i)$-entry of $B^*_0$ is equal to the
expression on the right in (\ref{eq:trace1}). This shows (\ref{eq:trace1}).
To obtain  (\ref{eq:trace2})--(\ref{eq:trace8}) we apply the $D_4$ action
to (\ref{eq:trace1}) and use Theorem  \ref{thm:rel}.
\end{proof}

\medskip

\begin{corollary}       \label{cor:nonzero}
In the equations (\ref{eq:trace1})--(\ref{eq:trace8}) 
each side is nonzero.
\end{corollary}

\begin{proof}
Immediate from Theorem \ref{thm:classify}(i).
\end{proof}

\medskip

We now present our trace formulae.

\medskip

\begin{theorem}      \label{thm:main}
For $1 \leq i \leq d$ we have
\begin{eqnarray}
 \varphi_i &=& (\theta^*_0 - \theta^*_i)
           \frac{\text{\rm tr}(\tau_i(A) E^*_0)}
                {\text{\rm tr}(\tau_{i-1}(A) E^*_0)},  \label{eq:1} \\
 \varphi_i &=& (\theta^*_d - \theta^*_{i-1})
           \frac{\text{\rm tr}(\eta_{d-i+1}(A) E^*_d)}
                {\text{\rm tr}(\eta_{d-i}(A) E^*_d)},   \label{eq:SD}\\
 \varphi_i &=& (\theta_0 - \theta_i)
           \frac{\text{\rm tr}(\tau^*_i(A^*) E_0)}
                {\text{\rm tr}(\tau^*_{i-1}(A^*) E_0)},  \label{eq:star}\\
 \varphi_i &=& (\theta_d - \theta_{i-1})
           \frac{\text{\rm tr}(\eta^*_{d-i+1}(A^*) E_d)}
                {\text{\rm tr}(\eta^*_{d-i}(A^*) E_d)},  \label{eq:SDstar}\\
 \phi_i &=& (\theta^*_0 - \theta^*_i)
           \frac{\text{\rm tr}(\eta_i(A) E^*_0)}
                {\text{\rm tr}(\eta_{i-1}(A) E^*_0)},  \label{eq:D}\\
 \phi_i &=& (\theta^*_d - \theta^*_{i-1})
           \frac{\text{\rm tr}(\tau_{d-i+1}(A) E^*_d)}
                {\text{\rm tr}(\tau_{d-i}(A) E^*_d)},  \label{eq:S}\\
 \phi_i &=& (\theta_0 - \theta_{d-i+1})
           \frac{\text{\rm tr}(\eta^*_{d-i+1}(A^*) E_0)}
                {\text{\rm tr}(\eta^*_{d-i}(A^*) E_0)},  \label{eq:Dstar}  \\
 \phi_i &=& (\theta_d - \theta_{d-i})
           \frac{\text{\rm tr}(\tau^*_i(A^*) E_d)}
                {\text{\rm tr}(\tau^*_{i-1}(A^*) E_d)}. \label{eq:Sstar}
\end{eqnarray}
We note that in (\ref{eq:1})--(\ref{eq:Sstar}) the denominator is nonzero
by Corollary \ref{cor:nonzero}.
\end{theorem}

\begin{proof}
Line (\ref{eq:1}) is routinely obtained from (\ref{eq:trace1}).
Lines (\ref{eq:SD})--(\ref{eq:Sstar}) are obtained from 
(\ref{eq:trace2})--(\ref{eq:trace8}) in a similar fashion.
\end{proof}

\medskip

\begin{remark}
The formula (\ref{eq:1}) was conjectured by the second author
\cite[Section 36]{T:survey}.
\end{remark}

\section{A transition formula}

Throughout this section 
$\Phi=(A; \{E_i\}_{i=0}^d; A^*; \{E^*_i\}_{i=0}^d)$ 
denotes a Leonard system in $\cal A$ and 
$p=(\{\theta_i\}_{i=0}^d; \{\theta^*_i\}_{i=0}^d;
        \{\varphi_i\}_{i=1}^d; \{\phi_i\}_{i=1}^d)$
denotes the corresponding parameter array.
We now define a parameter $q$.

\medskip

\begin{definition}
For $d \geq 3$ let $\beta$ denote the scalar in $\mathbb{K}$ such that 
$\beta+1$ is the common value of (\ref{eq:beta}).
For $d \leq 2$ let $\beta$ denote any scalar in $\mathbb{K}$.
Let $\overline{\mathbb{K}}$ denote the algebraic closure of $\mathbb{K}$.
Let $q$ denote a nonzero scalar in  $\overline{\mathbb{K}}$ such that
$\beta=q+q^{-1}$. 
In order to give our main idea clearly and to avoid limiting cases
we assume $q \neq 1$, $q \neq -1$.
\end{definition}

\medskip

\begin{lemma}  \cite[Lemma 9.2]{T:Leonard}     \label{lem:theta}
There exist scalars $\alpha$, $\mu$, $\nu$ in $\overline{\mathbb{K}}$
such that
\[
   \theta_i = \alpha + \mu q^{i} + \nu q^{-i}
   \qquad
   (0 \leq i \leq d).
\]
Replacing $q$ by $q^{-1}$ if necessary we may assume $\nu \neq 0$.
\end{lemma}

\medskip

\begin{lemma}        \label{lem:diff}
For $0 \leq i,j \leq d$ we have
\[
\theta_{i} - \theta_{j} = \nu q^{-i}(1-q^{i-j})(1- \mu \nu^{-1}q^{i+j}).
\]
\end{lemma}

\begin{proof}
Immediate from Lemma \ref{lem:theta}.
\end{proof}

\begin{corollary}     \label{cor:qinonzero}
For $1 \leq i \leq d$ we have $q^i \neq 1$.
\end{corollary}

\begin{proof}
By Lemma \ref{lem:diff} and since $\theta_0, \theta_1, \ldots, \theta_d$
are mutually distinct.
\end{proof}

\medskip

For a nonnegative integer $n$ we define
\begin{equation}    \label{eq:deffact}
  [n] = \frac{q^{n/2} - q^{-n/2}}{q^{1/2} - q^{-1/2}},
 \qquad
  [n]^! = [n][n-1] \cdots [2][1].
\end{equation}
We interpret $[0]^!=1$.
From Corollary \ref{cor:qinonzero} we find $[i] \neq 0$
for $1 \leq i \leq d$.
For nonnegative integers $r$, $s$, $t$ with $r+s+t \leq d$ we define
\begin{equation}     \label{eq:rst}
[r,s,t]= \frac{[r+s]^!\,[r+t]^!\,[s+t]^!}
              {[r]^!\,[s]^!\,[t]^!\,[r+s+t]^!}.
\end{equation}
Our goal for the rest of this section is to prove the following
theorem.

\medskip

\begin{theorem}          \label{thm:trans}
For $0 \leq i \leq d$ we have
\begin{equation}           \label{eq:trans}
  \eta_i = \sum_{h=0}^i [h,i-h,d-i]\eta_{i-h}(\theta_0) \tau_h.
\end{equation}
\end{theorem}

\medskip

We will use the following notation.
For $a \in \overline{\mathbb{K}}$ and for an integer $n \geq 0$ 
we define
\begin{equation}     \label{eq:defsymb}
   (a;q)_n = (1-a)(1-aq) \cdots (1-aq^{n-1}).
\end{equation}
We interpret $(a;q)_0=1$.

\medskip

\begin{lemma}       \label{lem:taueta}
For $0 \leq i,j \leq d$ we have
\begin{equation}        \label{eq:tauth}
 \tau_{i}(\theta_j)
   = (-1)^i \nu^i q^{-i(i-1)/2} (q^{-j};q)_i (\mu \nu^{-1} q^j; q)_i.
\end{equation}
Also
\begin{equation}         \label{eq:etath}
 \eta_{i}(\theta_j)
   = (-1)^i \mu^i q^{id-i(i-1)/2} (q^{j-d};q)_i (\mu^{-1}\nu q^{-j-d};q)_i
\end{equation}
if $\mu \neq 0$ and
\begin{equation}         \label{eq:etath2}
 \eta_{i}(\theta_j) = \nu^{i} q^{-ij} (q^{j-d};q)_i
\end{equation}
if $\mu=0$.
\end{lemma}

\begin{proof}
Routine verification 
using (\ref{eq:deftau}), (\ref{eq:defeta}),
(\ref{eq:defsymb}) and Lemma \ref{lem:diff}.
\end{proof}

\medskip

\begin{lemma}           \label{lem:fact}
For an integer $n \geq 0$ we have 
\[
  [n]^! = q^{-n(n-1)/4} (1-q)^{-n} (q;q)_n.
\]
\end{lemma}

\begin{proof}
Easily verified using  (\ref{eq:deffact}) and (\ref{eq:defsymb}).
\end{proof}

\medskip

We use the standard notation for $q$-hypergeometric series:
\[
_{r+1} \phi_r
 \left[
    \begin{matrix}
      a_1, a_2, \ldots, a_{r+1} \\
      b_1, \ldots, b_r
    \end{matrix}
    ; q, z
 \right]
 = \sum_{t=0}^{\infty} 
    \frac{(a_1;q)_t (a_2;q)_t \cdots (a_{r+1};q)_t z^t}
         {(b_1;q)_t \cdots (b_r;q)_t (q;q)_t}.
\]
The following is known as the $q$-Pfaff-Saalsch\"{u}tz identity.

\medskip

\begin{lemma}  \cite[p.~355]{GR}      \label{lem:Saal}
Let $a$, $b$, $c$ denote nonzero scalars in $\overline{\mathbb{K}}$
and let $n$ denote a nonnegative integer. Then
\[
_3 \phi_2
  \left[
    \begin{matrix}
      a, b, q^{-n} \\
      c, abc^{-1}q^{1-n}
    \end{matrix}
     ; q, q
  \right]
  = \frac{(ca^{-1};q)_n (cb^{-1};q)_n}
         {(c;q)_n (c a^{-1}b^{-1};q)_n}.
\]
\end{lemma}

\medskip\noindent
The following is known as the $q$-Chu-Vandermonde identity.

\medskip

\begin{lemma}  \cite[p.~354]{GR}      \label{lem:Van}
Let $a$, $c$ denote nonzero scalars in $\overline{\mathbb{K}}$
and let $n$ denote a nonnegative integer. Then
\[
_2 \phi_1
  \left[
    \begin{matrix}
      a, q^{-n} \\
      c
    \end{matrix}
     ; q, q
  \right]
  = \frac{(ca^{-1};q)_n}
         {(c;q)_n } \, a^n.
\]
\end{lemma}

\medskip

We are now ready to prove Theorem \ref{thm:trans}.

\medskip

\begin{proofof}{Thoerem \ref{thm:trans}}
Observe that each side of (\ref{eq:trans}) is a polynomial in $\lambda$
with degree at most $d$. By this and since 
$\theta_0, \theta_1, \ldots, \theta_d$ are mutually distinct,
it suffices to show (\ref{eq:trans}) holds for $\lambda = \theta_j$
$(0 \leq j \leq d)$. 
Let $j$ be given. We set $\lambda=\theta_j$ in (\ref{eq:trans}) and evaluate
the result using Lemmas \ref{lem:taueta}, \ref{lem:fact}. 
For $\mu \neq 0$ we simplify further using 
Lemma \ref{lem:Saal} (with $a=q^{-j}$, $b=\mu\nu^{-1}q^j$,
$c=q^{-d}$, $n=i$).
For $\mu = 0$ we simplify further using 
Lemma \ref{lem:Van} (with $a=q^{-j}$, $c=q^{-d}$, $n=i$).
In either case we routinely verify (26) holds at $\lambda = \theta_j$.
\end{proofof}

\section{Some formulae relating the first and  second split sequence}

In this section we prove a theorem that relates the first and second
split sequence.
In principle this theorem can be deduced directly from 
Theorem \ref{thm:classify};
however we present a more natural proof using our trace formula.

\medskip

\begin{theorem}      \label{thm:main2}
Let $\Phi=(A; \{E_i\}_{i=0}^d; A^*; \{E^*_i\}_{i=0}^d)$ 
denote a Leonard system in $\cal A$ and let
$p=(\{\theta_i\}_{i=0}^d; \{\theta^*_i\}_{i=0}^d;
        \{\varphi_i\}_{i=1}^d; \{\phi_i\}_{i=1}^d)$
denote the corresponding parameter array.
Then for $1 \leq i \leq d$ we have
\[
  \frac{\phi_{1} \phi_{2} \cdots \phi_{i}}
       {(\theta^*_{0}-\theta^*_{1})(\theta^*_{0}-\theta^*_{2}) \cdots
           (\theta^*_{0}-\theta^*_{i})}
= \sum_{h=0}^{i} 
 \frac{[h,i-h,d-i]\eta_{i-h}(\theta_0)  \varphi_{1} \varphi_{2} \cdots \varphi_{h}}
       {(\theta^*_{0}-\theta^*_{1})(\theta^*_{0}-\theta^*_{2}) \cdots
           (\theta^*_{0}-\theta^*_{h})},
\]
\[
  \frac{\varphi_{d}\varphi_{d-1}\cdots \varphi_{d-i+1}}
       {(\theta^*_{d}-\theta^*_{d-1})(\theta^*_{d}-\theta^*_{d-2}) \cdots
           (\theta^*_{d}-\theta^*_{d-i})}
= \sum_{h=0}^{i} 
  \frac{[h,i-h,d-i]\eta_{i-h}(\theta_0) \phi_{d} \phi_{d-1} \cdots \phi_{d-h+1}}
       {(\theta^*_{d}-\theta^*_{d-1})(\theta^*_{d}-\theta^*_{d-2}) \cdots
           (\theta^*_{d}-\theta^*_{d-h})},
\]
\[
  \frac{\varphi_{1} \varphi_{2} \cdots \varphi_{i}}
       {(\theta^*_{0}-\theta^*_{1})(\theta^*_{0}-\theta^*_{2}) \cdots
           (\theta^*_{0}-\theta^*_{i})}
= \sum_{h=0}^{i} 
  \frac{[h,i-h,d-i]\tau_{i-h}(\theta_d) \phi_{1} \phi_{2} \cdots \phi_{h}}
       {(\theta^*_{0}-\theta^*_{1})(\theta^*_{0}-\theta^*_{2}) \cdots
           (\theta^*_{0}-\theta^*_{h})},
\]
\[
  \frac{\phi_{d}\phi_{d-1} \cdots \phi_{d-i+1}}
       {(\theta^*_{d}-\theta^*_{d-1})(\theta^*_{d}-\theta^*_{d-2}) \cdots
           (\theta^*_{d}-\theta^*_{d-i})}
= \sum_{h=0}^{i}
  \frac{[h,i-h,d-i]\tau_{i-h}(\theta_d) \varphi_{d} \varphi_{d-1} \cdots \varphi_{d-h+1}}
       {(\theta^*_{d}-\theta^*_{d-1})(\theta^*_{d}-\theta^*_{d-2}) \cdots
           (\theta^*_{d}-\theta^*_{d-h})},
\]
\[
  \frac{\phi_{d} \phi_{d-1} \cdots \phi_{d-i+1}}
       {(\theta_{0}-\theta_{1})(\theta_{0}-\theta_{2}) \cdots
           (\theta_{0}-\theta_{i})}
= \sum_{h=0}^{i}
  \frac{[h,i-h,d-i]\eta^*_{i-h}(\theta^*_0)\varphi_{1} \varphi_{2} \cdots \varphi_{h}}
       {(\theta_{0}-\theta_{1})(\theta_{0}-\theta_{2}) \cdots
           (\theta_{0}-\theta_{h})},
\]
\[
  \frac{\varphi_{d} \varphi_{d-1} \cdots \varphi_{d-i+1}}
       {(\theta_{d}-\theta_{d-1})(\theta_{d}-\theta_{d-2}) \cdots
           (\theta_{d}-\theta_{d-i})}
= \sum_{h=0}^{i} 
  \frac{[h,i-h,d-i]\eta^*_{i-h}(\theta^*_0)\phi_{1} \phi_{2} \cdots \phi_{h}}
       {(\theta_{d}-\theta_{d-1})(\theta_{d}-\theta_{d-2}) \cdots
           (\theta_{d}-\theta_{d-h})},
\]
\[
  \frac{\varphi_{1} \varphi_{2} \cdots \varphi_{i}}
       {(\theta_{0}-\theta_{1})(\theta_{0}-\theta_{2}) \cdots
           (\theta_{0}-\theta_{i})}
= \sum_{h=0}^{i} 
  \frac{[h,i-h,d-i]\tau^*_{i-h}(\theta^*_d)\phi_{d} \phi_{d-1} \cdots \phi_{d-h+1}}
       {(\theta_{0}-\theta_{1})(\theta_{0}-\theta_{2}) \cdots
           (\theta_{0}-\theta_{h})},
\]
\[
  \frac{\phi_{1} \phi_{2} \cdots \phi_{i}}
       {(\theta_{d}-\theta_{d-1})(\theta_{d}-\theta_{d-2}) \cdots
           (\theta_{d}-\theta_{d-i})}
= \sum_{h=0}^{i}
  \frac{[h,i-h,d-i]\tau^*_{i-h}(\theta^*_d)\varphi_{d} \varphi_{d-1} \cdots \varphi_{d-h+1}}
       {(\theta_{d}-\theta_{d-1})(\theta_{d}-\theta_{d-2}) \cdots
           (\theta_{d}-\theta_{d-h})}.
\]
We are using the notation (\ref{eq:rst}).
\end{theorem}

\begin{proof}
Evaluating (\ref{eq:trans}) at $\lambda=A$ we find
\[
   \eta_i(A) = \sum_{h=0}^{i} [h,i-h,d-i]\eta_{i-h}(\theta_0) \tau_{h}(A).
\]
In the above line we multiply each side on the right by $E^*_0$ to obtain
\begin{equation}     \label{eq:aux}
   \eta_i(A) E^*_0 = \sum_{h=0}^{i} [h,i-h,d-i]\eta_{i-h}(\theta_0) \tau_{h}(A) E^*_0.
\end{equation}
In equation (\ref{eq:aux}) we take the trace of both sides
and use (\ref{eq:trace1}), (\ref{eq:trace3}) to obtain the first formula. 
Applying the $D_4$ action to the first formula
using Theorem \ref{thm:rel} we obtain the remaining formulae.
\end{proof}

\bibliographystyle{plain}

\bigskip\bigskip\noindent
Kazumasa Nomura\\
College of Liberal Arts and Sciences\\
Tokyo Medical and Dental University\\
Kohnodai, Ichikawa, 272-0827 Japan\\
email: nomura.las@tmd.ac.jp

\bigskip\noindent
Paul Terwilliger\\
Department of Mathematics\\
University of Wisconsin\\
480 Lincoln drive, Madison, Wisconsin, 53706 USA\\
email: terwilli@math.wisc.edu

\bigskip\noindent
{\bf Keywords.}
Leonard pair, tridiagonal pair, $q$-Racah polynomial, hypergeometric series.

\noindent
{\bf 2000 Mathematics Subject Classification}.
05E30, 05E35, 33C45, 33D45.

\end{document}